\documentclass[a4paper,11pt]{amsart}
\usepackage{amsmath}
\usepackage{amssymb}
\usepackage{amsthm}
\usepackage[all]{xy}
\usepackage[latin1]{inputenc}        
\usepackage[dvips]{graphics}
\usepackage[dvips]{graphicx}
\usepackage{amscd}
\usepackage{mathrsfs}
\usepackage[mathcal]{eucal}
\usepackage{fullpage}
\author{Giovanna Carnovale, Francesco Polizzi}
\address{Dipartimento di Matematica Pura ed Applicata, Università di Padova,
Via Trieste 63, 35121 Padova, Italy.} \email{carnoval@math.unipd.it}
\address{Dipartimento di Matematica, Università della Calabria, Via Pietro Bucci,
87036 Arcavacata di Rende (CS), Italy.}
\email{polizzi@mat.unical.it}
\title[Surfaces with $p_g=q=1$ isogenous to a product]
{The classification of surfaces with $p_g=q=1$
isogenous to a product of curves}
\date{\today}

\newtheorem{inizio}{Lemma}[section]

\newtheorem{proposition}[inizio]{Proposition}
\newtheorem{lemma}[inizio]{Lemma}

\newtheorem{definition}[inizio]{Definition}
\newtheorem*{teo}{Main Theorem}

\newtheorem{comp}[inizio]{Computational Fact}
\theoremstyle{definition}
\newtheorem{remark}[inizio]{Remark}

\newcommand{\lr}{\longrightarrow}

\newcommand{\mO}{\mathcal{O}}
\newcommand{\mZ}{\mathbb{Z}}

\newcommand{\Sn}{\sum_{j=1}^s \left(1- \frac{1}{ \;n_j} \right)}

\begin{document}

\subjclass[2000]{14J29 (primary), 14L30, 14Q99, 20F05} \keywords{Surfaces of general
type, isotrivial fibrations, actions of finite groups}

\abstract
A smooth, projective surface $S$ is said to be \emph{isogenous to a product}
if there exist two smooth curves $C$, $F$ and a finite group $G$ acting
 freely on $C \times F$
so that $S=(C \times F)/G$. In this paper we classify all surfaces with $p_g=q=1$ which are isogenous to a product.
\endabstract

\maketitle

\section{Introduction}
The classification of smooth, complex surfaces $S$ of general type with small birational
invariants is quite a natural problem in the framework of algebraic geometry. For
instance, one may want to understand the case where the Euler characteristic
$\chi(\mO_S)$ is $1$, that is, when the geometric genus $p_g(S)$ is equal to the
irregularity $q(S)$. All surfaces of general type with these invariants satisfy $p_g \leq
4$. In addition, if $p_g=q=4$ then the self-intersection $K_S^2$ of the canonical class
of $S$  is equal to $8$ and $S$ is the product of two genus $2$ curves, whereas if
$p_g=q=3$ then $K_S^2=6$ or $8$ and both cases are completely described (\cite{CCML98},
\cite{HP02}, \cite{Pir02}). On the other hand,
 surfaces of general type with $p_g=q=0, 1, 2$ are still
  far from being classified. We refer
  the reader to the survey paper \cite{BaCaPi06} for a recent account
  on this topic and a comprehensive list of references. \\
A natural way of producing interesting examples of algebraic surfaces is to construct
them as quotients of known ones by the action of a finite group. For instance
 Godeaux constructed in \cite{Go31} the first example of surface of general type
with vanishing geometric genus taking the quotient of a general quintic surface of
$\mathbb{P}^3$ by a free action of $\mathbb{Z}_5$. In line with this Beauville proposed
in \cite[p. 118]{Be2}  the construction of a surface of general type with $p_g=q=0$,
$K_S^2=8$
 as the quotient of a product of two curves $C$ and $F$ by the free action
 of a finite group $G$ whose
order is related to the genera $g(C)$ and $g(F)$ by the equality $|G|=(g(C)-1)(g(F)-1)$.
 Generalizing Beauville's example we say that a
 surface $S$ is  \emph{isogenous to a product} if $S=(C \times F )/G$, for
 $C$ and $F$ smooth curves and $G$ a finite group acting freely on $C \times F$. A
systematic study of these surfaces has been carried out in
\cite{Ca00}. They are of general type if and only if both $g(C)$ and
$g(F)$ are greater than or equal to $2$ and in this case $S$
admits a unique minimal realization where
 they are as small as possible. From now on, we tacitly assume that such a realization is
chosen, so that the genera of the curves and the group $G$ are invariants of $S$. The
action of $G$ can be seen to respect
 the product structure on $C \times F$. This means that such actions
 fall in two cases:
 the \emph{mixed} one, where there exists some element in $G$
 exchanging the two factors (in this situation $C$ and $F$ must be isomorphic)
 and the \emph{unmixed}
  one, where  $G$ acts faithfully on both
  $C$ and $F$ and diagonally on their product. \\
After \cite{Be2}, examples of surfaces isogenous to a product with $p_g=q=0$ appeared in
\cite{Par03} and  \cite{BaCa03}, and their complete classification was obtained in \cite{BaCaGr06}.\\
The next natural step is therefore the analysis of the case $p_g=q=1$. Surfaces of
general type with these invariants are the irregular ones with the lowest geometric genus
and for this reason it would be
  important to provide their complete description. So
  far, this has been obtained only in the cases $K_S^2=2,3$ (\cite{Ca},
  \cite{CaCi91}, \cite{CaCi93}, \cite{Pol05}, \cite{CaPi06}). \\
The goal of the present paper is to give the full list of surfaces
with $p_g=q=1$ that are isogenous to a product. Our work has to be
seen as the sequel to the article \cite{Pol07}, which describe all
unmixed cases with $G$ abelian and some unmixed examples with $G$
nonabelian. Apart from the complete list of the genera and groups
occurring, our paper contains the first examples of surfaces of
mixed type with $q=1$. The mixed cases turn out to be much less
frequent than the unmixed ones and, as when $p_g=q=0$, they occur
for only one value of the order of $G$.
 However, in contrast with what happens when $p_g=q=0$, the mixed cases
 do not correspond to the maximum
value of $|G|$ but appear for a rather small order, namely $|G|=16$. \\
Our classification procedure involves arguments from both geometry
and computational group theory. We will give here a brief account on
how the result is achieved.\\
If $S$ is any surface isogenous to a product and satisfying $p_g=q$ then $|G|$, $g(C)$,
$g(F)$ are related as in Beauville's example and we have $K^2_S=8$. Besides, if $p_g=q=1$
such surfaces are necessarily minimal and of general type (Lemma \ref{min-gen}).\\ If
$S=(C \times F)/G$ is of unmixed type then the two projections $\pi_C \colon C \times F
\lr C$, $\pi_F \colon C \times F \lr F$ induce two morphisms $\alpha \colon S \lr C/G$,
 $\beta \colon S \lr F/G$, whose smooth fibres are isomorphic to $F$ and $C$,
 respectively. Moreover, the geometry of $S$ is encoded in the geometry of the two coverings
 $h \colon C \lr C/G$, $f \colon F \lr F/G$ and the invariants of $S$
 impose strong restrictions on $g(C)$, $g(F)$
and $|G|$. Indeed we have $1=q(S)=g(C/G)+g(F/G)$ so we may assume that $E:=C/G$ is an
elliptic curve and $F/G \cong \mathbb{P}^1$. Then $\alpha \colon S \lr E$ is the Albanese
morphism of $S$ and the genus  $g_{\textrm{alb}}$ of the general Albanese fibre equals
$g(F)$. It is proven in \cite[Proposition 2.3]{Pol07} that $3 \leq
g(F) \leq 5$; in particular this allows us to control $|G|$.
The covers $f$ and $h$
are determined by two suitable systems of generators for $G$, that we call $\mathcal{V}$
and $\mathcal{W}$, respectively.
Besides, in order to obtain a free action of $G$ on $C \times F$
and a quotient $S$ with the desired invariants, $\mathcal{V}$ and $\mathcal{W}$ are
subject to strict conditions of combinatorial nature (Proposition
\ref{structureresult2}). The geometry imposes also strong
restrictions on  the possible $\mathcal{W}$ and the genus of
$C$, so the existence of  $\mathcal{V}$ and $\mathcal{W}$ and the
compatibility
conditions can be verified through a computer search.
It is worth mentioning that the classification of finite groups of
automorphisms acting on curves of genus lesser than or equal to
 $5$ could have also been retrieved
from the existing literature (\cite{Br90}, \cite{Ki03}, \cite{KuKi90},
\cite{KuKu90}).\\ \\
If $S=(C \times C)/G$ is of mixed type then the index two subgroup $G^\circ$ of $G$
 corresponding to transformations that do not exchange the
 coordinates in $C \times C$ acts faithfully on $C$. The quotient
  $E=C/G^\circ$ is isomorphic to the Albanese variety of $S$ and $g_{\textrm{alb}}=g(C)$
  (Proposition \ref{alb-mix}).
Moreover $g(C)$ may
  only be $5$, $7$ or $9$, hence $|G|$ is at most $64$ (Proposition
 \ref{class-mix}). The cover $h \colon C \lr E$ is determined by
  a suitable system of generators $\mathcal{V}$ for $G^\circ$ and since the action
  of $G$ on $C \times C$ is required to be free, combinatorial
  restrictions involving the elements of $\mathcal{V}$ and those of $G
  \setminus G^\circ$ have to be imposed (Proposition
  \ref{building-mixed}). Our
  classification is obtained by first listing those groups
  $G^\circ$ for which $\mathcal{V}$ exists and
  then by looking at the admissible extensions $G$ of $G^\circ$. We find that the only
  possibility occurring is for $g(C)=5$ so that $|G|$ is necessarily $16$
  (Propositions \ref{mixed-5}, \ref{no-mixed-7}, \ref{no-mixed-9}).\\ \\
In the last part of the paper we examine the
structure of the subset of the moduli space corresponding to surfaces
isogenous to a product with $p_g=q=1$.
It can be explicitly described by
calculating the number of orbits of the direct product of certain
mapping class groups with $\textrm{Aut}(G)$ acting on the set (of pairs) of
systems of generators (Proposition \ref{prop:moduli spaces}). In particular
it is possible to determine the number of irreducible connected components and
their respective dimensions, see the forthcoming article \cite{Pe08}. \\ \\
Our computations were carried out by using the computer algebra
program \verb|GAP4|, whose database includes all groups of order
less than $2000$, with the exception of $1024$ (see \cite{GAP4}).
For the reader's convenience we included the scripts in the Appendix. \\ \\
Now let us state the main result of this paper.
\begin{teo}
Let $S=(C \times F)/G$ be a surface with $p_g=q=1$, isogenous to a
product of curves. Then $S$ is minimal of general type and the occurrences
 for $g(F)$, $g(C)$, $G$, the dimension $D$ of the moduli space and the number $N$ of
 its connected components are precisely those in the table below.
%
%
\begin{table}[ht!]
\begin{center}
\begin{tabular}{|c|c|c|c|c|c|c|}
\hline
&&$ $ & \verb|IdSmall| & $ $ &&\\
$g(F)=g_{\emph{alb}}$&$g(C)$&$G$ & \verb|Group|$(G)$ & {\rm{Type}}&$D$&$N$ \\
\hline \hline
$3$&$3$&$(\mZ_2)^2$ &$G(4,2)$& {\rm unmixed} $(*)$ &$5$&$1$\\
\hline
$3$&$5$&$(\mZ_2)^3$&$G(8,5)$& {\rm unmixed} $(*)$&$4$&$1$\\
\hline
$3$&$5$&$\mZ_2\times \mZ_4$&$G(8,2)$& {\rm unmixed} $(*)$&$3$&$2$\\
\hline
$3$&$9$&$\mZ_2\times\mZ_8$&$G(16,5)$& {\rm unmixed} $(*)$&$2$&$1$\\
\hline
$3$& $5$& $D_4$ & $G(8,3)$ & {\rm unmixed}& $3$ &$1$\\
\hline
$3$& $7$& $D_6$ & $G(12,4)$ & {\rm unmixed} $(**)$ & $3$ &$1$\\
\hline $3$& $9$& $\mathbb{Z}_2 \times D_{4}$  & $G(16,11)$ & {\rm unmixed}& $3$ &$1$\\
\hline $3$& $13$&  $D_{2,12,5}$ & $G(24,5)$ & {\rm unmixed} & $2$& $1$\\
\hline $3$& $13$&  $\mZ_2 \times A_4$ & $G(24,13)$ & {\rm unmixed}& $2$ & $1$\\
\hline $3$& $13$&  $S_4$ & $G(24,12)$ & {\rm unmixed} & $2$ & $1$\\
\hline $3$& $17$& $\mZ_2 \ltimes (\mZ_2 \times \mZ_8)$ & $G(32,9)$ & {\rm unmixed} & $2$& $1$\\
\hline $3$& $25$& $\mZ_2 \times S_4$ & $G(48,48)$ &  {\rm unmixed} & $2$ & $1$\\
\hline $4$& $3$ & $S_3$ & $G(6,1)$ & {\rm unmixed} $(**)$ & $4$ & $1$\\
\hline $4$& $5$ & $D_6$ & $G(12,4)$ & {\rm unmixed} & $3$ & $1$\\
\hline $4$& $7$& $\mZ_3 \times S_3$ & $G(18,3)$ & {\rm unmixed} & $2$& $2$\\
\hline $4$& $7$& $\mZ_3 \times S_3$ & $G(18,3)$ & {\rm unmixed} & $1$ & $1$\\
\hline $4$& $9$& $S_4$ & $G(24,12)$ & {\rm unmixed}  $(**)$& $2$ &$1$\\
\hline $4$& $13$& $S_3 \times S_3$ & $G(36,10)$ & {\rm unmixed} & $1$ & $1$\\
\hline $4$& $13$& $\mZ_6 \times S_3$ & $G(36,12)$ & {\rm unmixed} & $1$ & $1$\\
\hline $4$& $13$& $\mZ_4 \ltimes (\mZ_3)^2$ & $G(36,9)$ & {\rm unmixed}& $1$ & $2$\\
\hline $4$& $21$& $A_5$ & $G(60,5)$ & {\rm unmixed} $(**)$ &$1$& $1$\\
\hline $4$& $25$& $\mZ_3\times S_4$ &$G(72,42)$ & {\rm unmixed}& $1$ & $1$\\
\hline $4$& $41$ & $S_5$ & $G(120,34)$ & {\rm unmixed} & $1$ & $1$\\
\hline $5$& $3$ & $D_4$ & $G(8,3)$ & {\rm unmixed} $(**)$ & $4$ & $1$\\
\hline $5$& $4$ & $A_4$ & $G(12,3)$ & {\rm unmixed} $(**)$ &$2$& $2$\\
\hline $5$& $5$ & $\mZ_4\ltimes(\mZ_2)^2$ & $G(16,3)$ & {\rm unmixed} & $2$ & $3$\\
\hline $5$& $7$ &$\mZ_2 \times A_4$ & $G(24,13)$ & {\rm unmixed} & $2$ & $2$ \\
\hline $5$& $7$ & $\mZ_2 \times A_4$ & $G(24,13)$ & {\rm unmixed} & $1$ & $1$\\
\hline $5$& $9$ & $\mZ_8 \ltimes (\mZ_2)^2$ & $G(32,5)$ & {\rm unmixed} & $1$ & $1$\\
\hline $5$& $9$ & $\mZ_2\ltimes D_{2,8,5}$ & $G(32,7)$ & {\rm unmixed} & $1$ & $1$\\
\hline $5$& $9$ & $\mZ_4\ltimes(\mZ_4\times\mZ_2)$ & $G(32,2)$ & {\rm unmixed} &$1$& $1$\\
\hline $5$& $9$ & $\mZ_4\ltimes(\mZ_2)^3$ & $G(32,6)$ & {\rm unmixed}& $1$ & $1$\\
\hline $5$& $13$& $(\mZ_2)^2\times A_4$ & $G(48,49)$ & {\rm unmixed} & $1$ & $1$\\
\hline $5$& $17$& $\mZ_4\ltimes(\mZ_2)^4$ & $G(64,32)$ & {\rm unmixed} & $1$ & $2$\\
\hline $5$& $21$& $\mZ_5\ltimes(\mZ_2)^4$ & $G(80,49)$ & {\rm unmixed} & $1$ & $2$\\
\hline $5$& $5$& $D_{2,8,3}$&$G(16,8)$ & {\rm mixed} & $2$ & $1$\\
\hline $5$& $5$& $D_{2,8,5}$&$G(16,6)$ & {\rm mixed} & $2$ & $3$\\
\hline $5$ & $5$ & $\mZ_4\ltimes(\mZ_2)^2$ & $G(16,3)$ & {\rm mixed} & $2$ & $1$\\
\hline
\end{tabular}
\end{center}
\end{table}
\end{teo}
Here \verb|IdSmallGroup|$(G)$ denotes the label of the group $G$ in the
\verb|GAP4| database of
small groups. The calculation of $N$ is due to Penegini and Rollenske,
see \cite{Pe08}, except for the cases marked with $(*)$, which were already
studied in \cite{Pol07}. The cases marked with $(**)$ also appeared
in \cite{Pol07}, but the computation of $N$ was missing. \\ \\
This work is organized as follows. \\
In Section \ref{sup-isogene} we collect the basic facts about surfaces isogenous
to a product, following the treatment given by Catanese in \cite{Ca00}
and we fix the algebraic setup. \\
In Section \ref{building} we apply the structure theorems of Catanese
to the case $p_g=q=1$ and this leads to Propositions \ref{structureresult2}
 and \ref{building-mixed}, that provide the translation of our classification
 problem from geometry to algebra. All these results are used in Sections \ref{classification-unm} and
\ref{classification-mix},
 which are the core of the paper and give the complete lists of the
 occurring groups and genera
 in the unmixed and mixed cases, respectively. \\
Finally, Section \ref{moduli spaces} is devoted to the description
of the moduli spaces.
\bigskip

$\mathbf{Notations \; and \; conventions}$. All varieties,
morphisms, etc. in this article are defined over $\mathbb{C}$. By
``surface'' we mean a projective, non-singular surface $S$, and for such
a surface $K_S$ denotes the canonical class, $p_g(S)=h^0(S, \; K_S)$ is the
\emph{geometric genus}, $q(S)=h^1(S, \; K_S)$ is the
\emph{irregularity} and $\chi(\mathcal{O}_S)=1-q(S)+p_g(S)$ is the
\emph{Euler characteristic}. Throughout the paper we use the
following notation for groups:
\begin{itemize}
\item $\mZ_n$: cyclic group of order $n$.
\item $D_{p,q,r}=\mathbb{Z}_p \ltimes \mathbb{Z}_q= \langle x,y \; |
\; x^p=y^q=1, \; xyx^{-1}=y^r \rangle$: split metacyclic group of order
$pq$. The group
$D_{2,n,-1}$ is the dihedral group of order $2n$ and it will be
denoted by $D_n$.
\item $S_n, \;A_n$: symmetric, alternating group on $n$ symbols.
\item If $x,\,y \in G$,
their commutator is defined as $[x,y]=xyx^{-1}y^{-1}$.
\item If $x \in G$ we denote by $\textrm{Int}_x$ the inner automorphism
of $G$ defined as $\textrm{Int}_x(g)=xgx^{-1}$.
\item \verb|IdSmallGroup|$(G)$ indicates the label of the group $G$ in
the  \verb|GAP4| database of small groups. For instance
\verb|IdSmallGroup|$(D_4)=G(8,3)$ and this means that $D_4$ is the
third in the
 list of groups of order $8$.
\end{itemize}
\bigskip

$\mathbf{Acknowledgements}.$ The authors wish to thank M. Penegini
and S. Rollenske for giving them a preliminary version of
\cite{Pe08} and for kindly allowing them to include their results in
the Main Theorem. Moreover they are indebted with the referee
for several valuable comments and suggestions to improve this article.

\section{Basic on surfaces isogenous to a product} \label{sup-isogene}

In this section we collect for the reader's convenience some
basic results on groups acting on curves and
surfaces isogenous to a product, referring to \cite{Ca00} for further details.

\begin{definition} \label{isog-prod}
A complex surface $S$ of general type is said to be \emph{isogenous to a product}
 if there exist two smooth curves $C$, $F$ and a finite group $G$ acting freely on $C \times F$
 so that $S=(C \times F)/G$.
\end{definition}
There are two cases: the \emph{unmixed} one, where $G$
 acts diagonally, and the \emph{mixed} one, where there exist elements of $G$ exchanging
 the two factors $($and then $C$, $F$ are isomorphic$)$.

In both cases, since the action of $G$ on $C \times F$ is free, we have
\begin{equation} \label{Order-G}
\begin{split}
K_S^2&=\frac{K_{C \times F}^2}{|G|}=\frac{8(g(C)-1)(g(F)-1)}{|G|} \\
\chi(\mO_S)&=\frac{\chi(\mO_{C \times F})}{|G|}=\frac{(g(C)-1)(g(F)-1)}{|G|},
\end{split}
\end{equation}
hence $K_S^2=8 \chi(\mO_S)$. \\
Let $C$, $F$ be curves of genus $\geq 2$.
Then the inclusion $\textrm{Aut}(C \times F) \supset \textrm{Aut}(C) \times \textrm{Aut}(F)$
 is an equality if $C$ and $F$ are not isomorphic, whereas
$\textrm{Aut}(C \times C) = \mZ_2 \ltimes  (\textrm{Aut}(C) \times \textrm{Aut}(C))$, the
$\mathbb{Z}_2$ being generated by the involution exchanging the two coordinates. If $S=(C
\times F)/G$ is a surface isogenous to a product, we will always consider its unique
\emph{minimal realization}. This means that
\begin{itemize}
\item in the unmixed case, we have $G \subset \textrm{Aut}(C)$ and
$G \subset \textrm{Aut}(F)$ (i.e. $G$ acts faithfully on both $C$ and $F$);
\item in the mixed case, where $C \cong F$, we have $G^\circ \subset \textrm{Aut}(C)$,
 for $G^\circ:=G \cap (\textrm{Aut}(C) \times \textrm{Aut}(C))$.
\end{itemize}
(See \cite[Corollary 3.9 and Remark 3.10]{Ca00}). \\

\begin{definition} \label{generating vect}
Let $G$ be a finite group and let
$\mathfrak{g}' \geq 0$, and $m_r \geq m_{r-1} \geq \ldots \geq m_1
\geq 2$
be integers. A \emph{generating vector} for $G$ of type
$(\mathfrak{g}' \; | \; m_1, \ldots ,m_r)$ is a $(2
\mathfrak{g}'+r)$-ple of elements
\begin{equation*}
\mathcal{V}=\{g_1, \ldots, g_r; \; h_1, \ldots, h_{2\mathfrak{g}'}
\}
\end{equation*}
such that:
the set $\mathcal{V}$ generates $G$;\;
$|g_i|=m_i$ and
$g_1g_2\cdots g_r \Pi_{i=1}^{\mathfrak{g}'} [h_i,h_{i+\mathfrak{g}'}]=1$.
If such a $\mathcal{V}$ exists, then $G$ is said to be
$(\mathfrak{g}' \; | \; m_1, \ldots ,m_r)$-\emph{generated}.
\end{definition}
For convenience we make abbreviations such as $(4 \;| \; 2^3, 3^2)$ for $(4 \; | \;
2,2,2,3,3)$ when we write down the type of the generating vector
$\mathcal{V}$. \\
By Riemann's existence theorem a finite group $G$ acts as a group of automorphisms of some compact
Riemann surface $X$ of genus $\mathfrak{g}$ with quotient a Riemann surface $Y$
of genus $\mathfrak{g}'$ if and only if there exist integers  $m_r \geq m_{r-1} \geq
\ldots \geq m_1 \geq 2$ such that $G$ is $(\mathfrak{g}'\; |\; m_1,
\ldots, m_r)$-generated and $\mathfrak{g}$, $\mathfrak{g}'$, $|G|$ and
the $m_i$ are related by the Riemann-Hurwitz formula.
Moreover, if $\mathcal{V}=\{g_1, \ldots, g_r; \; h_1, \ldots, h_{2\mathfrak{g}'}
\}$ is a generating vector for $G$, the subgroups $\langle g_i\rangle$
and their conjugates are precisely the nontrivial stabilizers of the $G$-action
([Br90, Section 2], [Bre00, Chapter 3], \cite{H71}).
The description of surfaces isogenous to a product can be therefore reduced
to finding suitable generating vectors. Requiring that $S$
has given invariants $p_g$ and $q$ imposes numerical restrictions on
the order of the group $G$ and the genus of the curves $C$ and $F$.
Our goal is to classify all surfaces with $p_g=q=1$ isogenous to a product.
The aim of the next
section is to translate this classification problem
from geometry to algebra.

\section{The case $p_g=q=1$. Building data} \label{building}

\begin{lemma} \label{min-gen}
Let $S=(C \times F)/G$ be a surface isogenous to a product with $p_g=q=1$.
Then
\begin{itemize}
\item[$(i)$] $K_S^2=8$.
\item[$(ii)$] $|G|=(g(C)-1)(g(F)-1)$.
\item[$(iii)$] $S$ is a minimal surface of general type.
\end{itemize}
\end{lemma}
\begin{proof}
Claims $(i)$ and $(ii)$ follow from \eqref{Order-G}. Now let us
consider $(iii)$. Since $C \times F$ is minimal and the cover  $C
\times F \lr S$ is étale, $S$ is minimal as well. Moreover $(ii)$
implies either $g(C)=g(F)=0$ or $g(C) \geq 2$, $g(F) \geq 2$. The
first case is impossible otherwise $S= \mathbb{P}^1 \times
\mathbb{P}^1$ and $p_g=q=0$; thus the second case occurs, hence $S$
is of general type.
\end{proof}

\subsection{Unmixed case}\label{unmixed-structure}

If $S=(C \times F) /G$ is a surface with $p_g=q=1$,
isogenous to an unmixed product, then $g(C) \geq 3$, $g(F) \geq 3$
and up to exchanging $F$ and $C$ one may assume
 $F/G \cong \mathbb{P}^1$ and $C/G \cong E$, where $E$ is
 an elliptic curve. Moreover $\alpha \colon S  \lr C/G$ is the Albanese
 morphism of $S$ and
$g_{\emph{alb}}=g(F)$, see \cite[Proposition 2.2]{Pol07}.
This leads to
\begin{proposition}$($\cite[Proposition 3.1]{Pol07}$)$ \label{structureresult2}
Let $G$ be a finite group which is both $(0\; | \; m_1, \ldots
,m_r)$ and $(1\; | \; n_1, \ldots, n_s)$-generated, with generating
vectors $\mathcal{V}=\{g_1, \ldots, g_r \}$ and
$\mathcal{W}=\{\ell_1, \ldots, \ell_s; \; h_1,h_2 \}$, respectively.
Let $g(F), \; g(C)$ be the positive integers defined by the Riemann-Hurwitz relations
\begin{equation}\label{generi1}
\begin{split}
2g(F)-2
 &=|G|\bigg(-2+\sum_{i=1}^r \bigg( 1- \frac{1}{\;m_i} \bigg)
\bigg),\qquad
\\
2g(C)-2 & =|G| \sum_{j=1}^s \bigg(1-\frac{1}{\;n_j} \bigg).
\end{split}
\end{equation}
Assume moreover that
$g(C) \geq 3$, $g(F) \geq3$,
$|G|=(g(C)-1)(g(F)-1)$ and
\begin{equation*}(U) \quad \quad \quad \left(\bigcup_{\sigma\in
    G}\bigcup_{i=1}^r \langle \sigma g_i \sigma^{-1} \rangle \right) \cap \left(\bigcup_{\sigma\in
    G}\bigcup_{j=1}^s \langle  \sigma \ell_j \sigma^{-1} \rangle \right)=\{1_G \}.
\end{equation*}
Then there is a free, diagonal action of $G$ on $C \times F$ such that the quotient $S=(C
\times F)/G$ is a minimal surface of general type with $p_g=q=1$, $K_S^2=8$. Conversely,
 every surface with $p_g=q=1$, isogenous to an unmixed product, arises in this way.
\end{proposition}
Here, condition $(U)$ ensures that the $G$-action on $C\times F$ is
free. \\
Set $\mathbf{m}:=(m_1, \ldots, m_r)$ and $\mathbf{n}:=(n_1, \ldots, n_s)$;
if $S=(C \times F)/G$ is a surface with $p_g=q=1$ which is
constructed by using the
recipe in Proposition \ref{structureresult2}, it will be called an
\emph{unmixed surface of type} $(G, \, \mathbf{m}, \, \mathbf{n})$.
\begin{proposition} $($\cite[Proposition 2.3]{Pol07}$)$ \label{various-cases}
Let $S=(C \times F)/G$ be an unmixed surface of type $(G, \,
\mathbf{m}, \,\mathbf{n})$. Then there are exactly the following
possibilities:
\begin{itemize}
\item[$(a)$] $g(F)=3, \quad \mathbf{n}=(2^2)$
\item[$(b)$] $g(F)=4, \quad \mathbf{n}=(3)$
\item[$(c)$] $g(F)=5, \quad \mathbf{n}=(2)$.
\end{itemize}
\end{proposition}
The following lemma gives a restriction on $\mathbf{m}$ instead.
\begin{lemma}\label{pollicino}Let $S=(C \times F)/G$ be an unmixed surface of
type $(G, \, \mathbf{m}, \, \mathbf{n})$. Then every $m_i$ divides
  $\frac{|G|}{(g(F)-1)}$.
\end{lemma}
\begin{proof}
Since $\langle g_i \rangle$ is a stabilizer for the $G$-action on $F$
and since $G$ acts freely on $(C \times F)$, the subgroup
 $\langle g_i\rangle\cong\mZ_{m_i}$ acts
  freely on $C$.  By Riemann-Hurwitz formula applied to the cover $C \lr
  C/\langle g_i\rangle$ we have
$g(C)-1=m_i( g(C/\langle g_i \rangle)-1)$. Thus $m_i$ divides
$g(C)-1=\frac{|G|}{(g(F)-1)}$.
\end{proof}

\subsection{Mixed case}

\begin{proposition} \label{alb-mix}
Let $S=(C \times C)/G$ be a surface with $p_g=q=1$ isogenous to a mixed
product. Then $E:=C/G^\circ$ is an elliptic curve isomorphic to the Albanese variety of
$S$.
\end{proposition}
\begin{proof}
We have (see \cite[Proposition 3.15]{Ca00})
\begin{equation*}
\begin{split}
\mathbb{C}=H^0(\Omega^1_S) & =(H^0(\Omega^1_C) \oplus H^0(\Omega^1_C))^G
  = (H^0(\Omega^1_C)^{G^\circ} \oplus H^0(\Omega^1_C)^{G^\circ})^{G/G^{\circ}} \\
& = (H^0(\Omega^1_E) \oplus H^0(\Omega^1_E))^{G/G^{\circ}}.
\end{split}
\end{equation*}
Since $S$ is of mixed type, the quotient $\mathbb{Z}_2=G/G^\circ$
 exchanges the last two summands, whence $h^0(\Omega_E^1)=1$.
Thus $E$ is an elliptic curve and there is a commutative diagram
\begin{equation} \label{diagram-alb}
\xymatrix{ C \times C \ar[r]^{\rho}  \ar[d]^{\pi} &
 E \times E \ar[d]^{\varepsilon} \\
 S \ar[r]^{\hat{\rho}} \ar[dr]^{\alpha} &
 E^{(2)} \ar[d]^{\hat{\alpha}} \\
  & E}
\end{equation}
showing that the Albanese morphism $\alpha$ of $S$ factors through the
Abel-Jacobi map $\hat{\alpha}$ of the double symmetric product $E^{(2)}$ of $E$.
\end{proof}
By Lemma \ref{min-gen} we have $|G|=(g(C)-1)^2$. In this case
\cite[Proposition 3.16]{Ca00} becomes
\begin{proposition} \label{building-mixed}
Assume that $G^\circ$ is a $(1\; | \; n_1, \ldots, n_s)$-generated
finite group
with generating vector
 $\mathcal{V}=\{\ell_1, \ldots, \ell_s; \; h_1, h_2 \}$
and that there is a nonsplit extension
\begin{equation} \label{nonsplit}
1 \lr G^\circ \lr G \lr \mZ_2 \lr 1
\end{equation}
which gives an involution $[\varphi]$ in \emph{Out}$(G^\circ)$.
Let $g(C)\in \mathbb{N}$ be defined by the Riemann-Hurwitz relation $2g(C)-2=|G^\circ| \Sn$.
Assume, in addition, that $|G|=(g(C)-1)^2$ and that
\begin{itemize}
\item[$(M1)$]
for all $g \in G \setminus G^\circ$ we have
\begin{equation*}
\{\ell_1, \ldots, \ell_s \} \cap \{ g \ell_1 g^{-1}, \ldots, g
\ell_s g^{-1} \}=\emptyset;
\end{equation*}
\item[$(M2)$] for all $g \in G \setminus G^\circ$ we have
\begin{equation*}
g^2 \notin \bigcup _{j=1}^s \bigcup_{\sigma \in G^\circ} \langle
\sigma \ell_j \sigma^{-1} \rangle.
\end{equation*}
\end{itemize}
Then there is a free, mixed action of $G$ on $C \times C$ such that the quotient $S=(C
\times C)/G$ is a minimal surface of
general type with $p_g=q=1$, $K_S^2=8$. \\
Conversely, every surface $S$ with $p_g=q=1$,
isogenous to a mixed product, arises in this way.
\end{proposition}
Here, conditions $(M1)$ and $(M2)$ ensure that the $G$-action on
$C\times C$ is free.

\begin{remark} \label{NO-ABELIAN}
 The surface $S$ is not covered by elliptic curves because it is of
 general type (Lemma \ref{min-gen}), so
 the map $C\lr C/G^\circ=E$ is ramified. Therefore condition
$(M1)$ implies that $G$ is not abelian.
\end{remark}

\begin{remark} \label{elements-order-2}
The exact sequence (\ref{nonsplit}) is non split if and only if the
number of elements of order $2$ in $G$ equals the number of elements
of order $2$ in $G^\circ$.
\end{remark}
\begin{proposition} \label{alb-mix-2}
Let $S=(C \times C)/G$ be a surface with $p_g=q=1$, isogenous to a mixed product.
Then $g_{\emph{alb}}=g(C)$.
\end{proposition}
\begin{proof}
Let us look at diagram \eqref{diagram-alb}. The Abel-Jacobi map $\hat{\alpha}$ gives to
$E^{(2)}$ the structure of a $\mathbb{P}^1$-bundle over $E$ (\cite{CaCi93}); let
$\mathfrak{f}$ be the generic fibre of this bundle and $F^*:=\rho^* \varepsilon^*
(\mathfrak{f})$. If $F_{\textrm{alb}}$ is the generic Albanese fibre of $S$ we have
$F_{\textrm{alb}}=\pi(F^*)$. Let $\mathbf{n}=(n_1,
\ldots, n_s)$ be such that $G^\circ$ is $(1 \,| \,n_1, \ldots n_s)$-generated and $2g(C)-2=|G^\circ| \Sn$.
 The $(G^\circ \times G^\circ)$-cover $\rho$ is branched exactly along the union of
 $s$ ``horizontal'' copies of $E$ and $s$ ``vertical'' copies of $E$; moreover
 for each $i$ there are one horizontal copy and one vertical copy
 whose branching number is $n_i$. Since $\varepsilon^*(\mathfrak{f})$
 is an elliptic curve that
  intersects all these copies of $E$ transversally in one point, by Riemann-Hurwitz
 formula applied to $F^* \lr \varepsilon^*(\mathfrak{f})$ we obtain
\begin{equation*}
2g(F^*)-2=|G^\circ|^2 \cdot \sum_{j=1}^s  2 \left(1- \frac{1}{ \;n_j} \right).
\end{equation*}
On the other hand the $G$-cover $\pi$ is étale, so we have
\begin{equation*}
\begin{split}
2g(F_{\textrm{alb}})-2= \frac{1}{|G|}(2g(F^*)-2)
& =|G^\circ| \Sn
\\
& =2g(C)-2,
\end{split}
\end{equation*}
whence $g_{\textrm{alb}}=g(C)$.
\end{proof}
If $S=(C \times C)/G$ is a surface with $p_g=q=1$ which is constructed
 by using the recipe of Proposition \ref{building-mixed},
it will be called a \emph{mixed surface of type} $(G, \,
\mathbf{n})$. The analogue of Proposition \ref{various-cases} in the
mixed case is
\begin{proposition} \label{class-mix}
Let $S=(C \times C)/G$ be a mixed surface of type $(G, \, \mathbf{n})$.
Then there are at most the following
possibilities:
\begin{itemize}
\item $g(C)=5, \quad \mathbf{n}=(2^2), \quad |G|=16$;
\item $g(C)=7, \quad \mathbf{n}=(3), \quad \;\,|G|=36$;
\item $g(C)=9, \quad \mathbf{n}=(2), \quad \; \,|G|=64$.
\end{itemize}
\end{proposition}
\begin{proof}
By Proposition \ref{building-mixed} we have $2g(C)-2= |G^\circ| \Sn$ and $|G^\circ|=\frac{1}{2}(g(C)-1)^2$,
so $g(C)$ must be odd and we obtain
$4=(g(C)-1) \Sn$.
Therefore $4 \geq \frac{1}{2}(g(C)-1)$ and the only possibilities are $g(C)=3,5,7,9$. \\
The case $g(C)=3$ is ruled out because $G$ cannot be abelian by
Remark \ref{NO-ABELIAN}.\\
If $g(C)=5$ then $\Sn=1$, so $\mathbf{n}=(2^2)$ and $|G|=16$. \\
If $g(C)=7$ then $\Sn=\frac{2}{3}$, so $\mathbf{n}=(3)$ and
$|G|=36$. \\
If $g(C)=9$ then $\Sn=\frac{1}{2}$, so $\mathbf{n}=(2)$ and
$|G|=64$.
\end{proof}
We will see in Section \ref{class-mix} that only the case $g(C)=5$ actually occurs.
%

\section{The unmixed case} \label{classification-unm}
The classification of surfaces of general type with $p_g=q=1$
 isogenous to an unmixed product is carried out in \cite{Pol07} when
the group $G$ is abelian. Therefore in this section we assume that $G$
 is nonabelian. \\

Following \cite[Section 1.2]{BaCaGr06}, for an $r$-ple
$\mathbf{m}=(m_1, \ldots, m_r) \in \mathbb{N}^r$ we set
\begin{equation*}
\Theta(\mathbf{m})
:=-2+ \sum_{i=1}^r \left(1- \frac{1}{ \;m_i} \right),
\quad\quad
\alpha(\mathbf{m})
:= \frac{2}{\Theta(\mathbf{m})}.
\end{equation*}
If $S$ is an unmixed surface of
type $(G, \, \mathbf{m}, \, \mathbf{n})$ then we necessarily have $2
\leq m_1 \leq \ldots \leq m_r$ and $\Theta(\mathbf{m}) >0$. Besides,
by Proposition \ref{structureresult2} we have $\alpha(\mathbf{m})=\frac{|G|}{g(F)-1}=g(C)-1 \in
\mathbb{N}$ and by Lemma \ref{pollicino} each integer $m_i$ divides
$\alpha(\mathbf{m})$. Then we get
%
\begin{proposition} \label{finite-tuples}
Let $S=(C \times F)/G$ be a surface with $p_g=q=1$ isogenous to an
unmixed product of type $(G, \, \mathbf{m}, \, \mathbf{n})$. Then
the possibilities for $\mathbf{m}$ and $\alpha(\mathbf{m})$, written
in the format $\mathbf{m}_{\alpha(\mathbf{m})}$, lie in the set
$\mathcal{T}$ below:
%
\begin{equation*}
{\mathcal T}=\left\{
\begin{matrix}
( 2 ,  3 , 7 )_{ 84 }, & ( 2 ,  3 , 8 )_{ 48 }, & ( 2 ,  4 , 5 )_{
40 }, & ( 2 ,  3 , 9 )_{ 36 }, & ( 2 ,  3 , 10 )_{ 30 }, & ( 2 ,
3 , 12 )_{ 24 }, \cr
( 2 ,  4 , 6 )_{ 24 }, & ( 3^2 , 4 )_{ 24 }, &
( 2 ,  5^2  )_{ 20 }, & ( 2 ,  3 , 18 )_{ 18 }, &  ( 2 ,  4 , 8
)_{ 16 }, & ( 3^2 , 5 )_{ 15 }, \cr
( 2 ,  4 , 12 )_{ 12 }, & ( 2 ,6^2 )_{ 12 }, & ( 3^2 , 6 )_{ 12
}, & ( 3 ,  4^2 )_{ 12 },
& ( 2 ,  5 , 10 )_{ 10 }, & ( 3^2 , 9 )_{ 9 }, \cr
 ( 2 ,  8^2 )_{8 }, & ( 4^3 )_{ 8 }, & ( 3 ,  6^2 )_{ 6 },  & ( 5^3)_{ 5 },  &( 2^3, 3 )_{ 12 },& ( 2^3 , 4 )_{ 8 },\cr
( 2^3 , 6)_{ 6 }, & ( 2^2 , 3^2 )_{ 6 }, & ( 2^2 , 4^2 )_{ 4 }, & (
3^4 )_{ 3 }, & ( 2^5 )_{4}, & ( 2^6)_{2} \end{matrix} \right\}
\end{equation*}
\end{proposition}
\begin{proof}
This follows combining \cite[Proposition 1.4]{BaCaGr06} with Lemma \ref{pollicino}.\end{proof}
By abuse of notation, we write $\mathbf{m}\in \mathcal{T}$
instead of $\mathbf{m}_{\alpha(\mathbf{m})}\in \mathcal{T}$.\\

Now we analyze the three cases in Proposition \ref{various-cases}
separately, according to the value of $g(F)$.
Note that if
$g(F)=3$, $4$, $5$ then $|\textrm{Aut}(F)| \leq 168$, $120$, $192$,
respectively (\cite[p. 91]{Bre00}).


\begin{proposition} \label{g(F)=3nonmixed}
If $g(F)=3$ we have precisely the following possibilities.
\begin{table}[ht!]
\begin{center}
\begin{tabular}{|c|c|c|}
\hline
$ $ & \verb|IdSmall| & $ $ \\
$G$ & \verb|Group|$(G)$ & $\mathbf{m}$ \\
\hline \hline
$D_4$ & $G(8,3)$ & $(2^2,4^2)$\\
\hline $D_6$ & $G(12,4)$ & $(2^3,6)$ \\
\hline $\mathbb{Z}_2 \times D_{4}$  & $G(16,11)$ & $(2^3,4)$
\\
\hline  $D_{2,12,5}$ & $G(24,5)$ & $(2,4,12)$ \\ \hline $\mZ_2
\times A_4$ & $G(24,13)$ & $(2,6^2)$ \\ \hline
 $S_4$ & $G(24,12)$ &
$(3,4^2)$
\\ \hline
$\mZ_2 \ltimes (\mZ_2 \times \mZ_8)$ & $G(32,9)$ & $(2,4,8)$ \\
\hline $\mZ_2 \times S_4$ & $G(48,48)$ & $(2,4,6)$ \\ \hline
\end{tabular}
\end{center}
\end{table}
\end{proposition}
\begin{proof}
Since $\mathbf{n}=(2^2)$ it follows that $G$ is $(1 \; | \;
2^2)$-generated and by the second relation in $(\ref{generi1})$ we
have $|G|=2(g(C)-1)$. So we must describe all unmixed surfaces of
type $(G, \mathbf{m}, \mathbf{n})$ with $\mathbf{m} \in
\mathcal{T}$, $\mathbf{n}=(2^2)$ and $|G|=2\alpha(\mathbf{m})$. By a
 computer search through the $r$-tuples in Proposition
\ref{finite-tuples} we can therefore list all possibilities, proving
our statement. See the \verb|GAP4| script $1$ in the Appendix to
see how this procedure applies to an explicit example.

\end{proof}


\begin{proposition} \label{genere4}
If $g(F)=4$ we have precisely the following possibilities.
\begin{table} [ht!]
\begin{center}
\begin{tabular}{|c|c|c|}
\hline
$ $ & \verb|IdSmall| & $ $ \\
$G$ & \verb|Group|$(G)$ & $\mathbf{m}$ \\
\hline \hline $S_3$ & $G(6,1)$ & $(2^6)$ \\ \hline $D_6$ & $G(12,4)$
& $(2^5)$ \\ \hline
 $\mZ_3 \times S_3$ & $G(18,3)$ & $(2^2,3^2)$ \\ \hline
 $\mZ_3 \times S_3$ & $G(18,3)$ & $(3,6^2)$ \\ \hline
$S_4$ & $G(24,12)$ & $(2^3,4)$ \\
\hline $S_3 \times S_3$ & $G(36,10)$ & $(2,6^2)$ \\
\hline $\mZ_6\times S_3$ & $G(36,12)$ & $(2,6^2)$ \\
\hline $\mZ_4 \ltimes (\mZ_3)^2$ & $G(36,9)$ & $(3,4^2)$ \\
\hline $A_5$ & $G(60,5)$ & $(2,5^2)$ \\
\hline $\mZ_3\times S_4$ &$G(72,42)$ & $(2,3,12)$
\\ \hline $S_5$ & $G(120,34)$ & $(2,4,5)$ \\ \hline
\end{tabular}
\end{center}
\end{table}
\end{proposition}
\begin{proof}
Since $\mathbf{n}=(3)$ it follows that $G$ is $(1\; | \;
3)$-generated and by the second relation in \eqref{generi1} we have
$|G|=3(g(C)-1)$. Therefore our statement can be proven searching by
computer calculation all unmixed surfaces of type $(G, \mathbf{m},
\mathbf{n})$ with $\mathbf{m} \in \mathcal{T}$, $\mathbf{n}=(3)$,
 $|G|=3\alpha(\mathbf{m})$ and $\alpha(\mathbf{m}) \leq 40$.
\end{proof}


\begin{proposition}
If $g(F)=5$ we have precisely the following possibilities.
\begin{table} [ht!]
\begin{center}
\begin{tabular}{|c|c|c|}
\hline
$ $ & \verb|IdSmall| & $ $ \\
$G$ & \verb|Group|$(G)$ & $\mathbf{m}$ \\
\hline \hline $D_4$ & $G(8,3)$ & $(2^6)$ \\
\hline $A_4$ & $G(12,3)$ & $(3^4)$ \\
\hline $\mZ_4\ltimes(\mZ_2)^2$ & $G(16,3)$ & $(2^2, 4^2)$ \\
\hline$\mZ_2 \times A_4$ & $G(24,13)$ & $(2^2, 3^2)$ \\
\hline  $\mZ_2 \times A_4$ & $G(24,13)$ & $(3,6^2)$ \\
\hline $\mZ_8 \ltimes (\mZ_2)^2$ & $G(32,5)$ & $(2,8^2)$ \\
\hline $\mZ_2\ltimes D_{2,8,5}$ & $G(32,7)$ & $(2,8^2)$ \\
\hline $\mZ_4\ltimes(\mZ_4\times\mZ_2)$ & $G(32,2)$ & $(4^3)$ \\
\hline $\mZ_4\ltimes(\mZ_2)^3$ & $G(32,6)$ & $(4^3)$ \\
\hline $(\mZ_2)^2\times A_4$ & $G(48,49)$ & $(2,6^2)$ \\
\hline $\mZ_4\ltimes(\mZ_2)^4$ & $G(64,32)$ & $(2,4,8)$ \\
\hline $\mZ_5\ltimes(\mZ_2)^4$ & $G(80,49)$ & $(2,5^2)$ \\
\hline
\end{tabular}
\end{center}
\end{table}

\end{proposition}
\begin{proof}
Since $\mathbf{n}=(2)$, it follows that $G$ is $(1\; | \;
2)$-generated and by the second relation in (\ref{generi1}) we have
$|G|=4(g(C)-1)$. Therefore our statement can be proven searching by
computer calculation all unmixed surfaces of type $(G, \mathbf{m},
\mathbf{n})$ with $\mathbf{m} \in \mathcal{T}$, $\mathbf{n}=(2)$,
$|G|=4\alpha(\mathbf{m})$ and $\alpha(\mathbf{m}) \leq 48$.
\end{proof}

\section{The mixed case} \label{classification-mix}
In this section we use Proposition \ref{building-mixed} in order to classify
the surfaces with $p_g=q=1$ isogenous to a mixed product. By Proposition \ref{class-mix}
 we have $g(C)=5$,  $7$ or $9$. Let us consider the three cases
 separately.

\newpage

\subsection{The case $g(C)=5, \; |G|=16$}
\begin{proposition} \label{mixed-5}
If $g(C)=5, \; |G|=16$ we have precisely the following possibilities.
\begin{table}[ht!]
\begin{center}
\begin{tabular}{|c|c|c|c|}
\hline
$ $ & \verb|IdSmall| & $ $ & \verb|IdSmall| \\
$G^{\circ}$ & \verb|Group|$(G^\circ)$ & $G$& \verb|Group|$(G)$ \\
\hline \hline
$D_4$&$G(8,3)$&$D_{2,8,3}$&$G(16,8)$\\
\hline
$\mZ_2 \times \mZ_4$&$G(8,2)$&$D_{2,8,5}$&$G(16,6)$\\
\hline
$(\mZ_2)^3$&$G(8,5)$&$\mZ_4 \ltimes (\mZ_2)^2$&$G(16,3)$\\
\hline
\end{tabular}
\end{center}
\end{table}
\end{proposition}
\begin{proof}
In this case $\mathbf{n}=(2^2)$, so our first task is to find all nonsplit sequences of
type \eqref{nonsplit} for which $G^\circ$ is a $(1\;|\;2^2)$-generated group of order
$8$. The three abelian groups of order $8$ and $D_4$ are $(1\;|\;2^2)$-generated whereas
the quaternion group $Q_8$ is not. \\
Since $\mZ_8$ has only one element $\ell$ of order $2$,
condition $(M1)$ in Proposition \ref{building-mixed} cannot be satisfied for any choice
of $\mathcal{V}$. By Remark \ref{NO-ABELIAN} we are left to analyze the possible
embeddings of $\mZ_2\times\mZ_4$, $D_4$ and $(\mZ_2)^3$ in nonabelian groups of order
$16$. The groups $\mZ_2 \times\mZ_4$, $D_4$ and $(\mZ_2)^3$ have $3$, $5$ and $7$
elements of order $2$, respectively. Therefore if $n_2$ denotes the number of elements of
order $2$ in $G$, by Remark \ref{elements-order-2} we must consider only those groups $G$
of order $16$ with $n_2\in\{3,\,5,\,7\}$. The nonabelian groups of order $16$ with
$n_2=3$ are $D_{2,8,5},\,\mZ_2\times Q_8$ and $D_{4,4,-1}$ and they all contain a copy of
$\mZ_2\times\mZ_4$. The only nonabelian group of order $16$ with $n_2=5$ is $D_{2,8,3}$
and it contains a subgroup isomorphic to $D_4$. The  nonabelian groups of order $16$ with
$n_2=7$ are $\mZ_4\ltimes(\mZ_2)^2=G(16,3)$ and $\mZ_2\ltimes Q_8$, and only the former
contains a subgroup isomorphic to
$(\mZ_2)^3$ (cfr. \cite{Wi05}). \\
Summarizing, we are left
with the following cases:
$$
\begin{array}{c|c}
G^\circ&G\\
\hline
D_4&D_{2,8,3}\\
\mZ_2\times\mZ_4&D_{2,8,5}\\
\mZ_2\times\mZ_4&\mZ_2\times Q_8\\
\mZ_2\times\mZ_4&D_{4,4,-1}\\
(\mZ_2)^3 &\mZ_4\ltimes(\mZ_2)^2\\
\end{array}
$$
Let us analyze them separately.\\ \\
$\bullet$ $G^\circ=D_4, \; \; G=D_{2,8,3}=\langle
x,y \; | \; x^2=y^8=1,\;xyx^{-1}=y^{3}\;\rangle$\\
We consider the subgroup $G^\circ:=\langle x,\,y^2\rangle\cong D_4$. Set
$\ell_1=\ell_2=x$ and $h_1=h_2=y^2$.
Condition $(M1)$ holds because $C_G(x)= \langle x, \, y^4 \rangle
\subset G^\circ$. Condition $(M2)$ is satisfied because the conjugacy
class of $x$ in $G^\circ$ is contained in the coset $x\langle
y^2\rangle$ while for every $g\in yG^\circ$ we have $g^2\in\langle
y\rangle$. Therefore this case occurs by Proposition \ref{building-mixed}.\\ \\
$\bullet$  $G^\circ=\mZ_2\times\mZ_4, \; \; G=D_{2,8,5}=\langle x,y\;
| \; x^2=y^8=1,
 \; xyx^{-1}=y^5 \rangle$\\
%
%
We consider the subgroup $G^\circ:=\langle x,\, y^2\rangle\cong\mZ_2\times\mZ_4$.
Set
$\ell_1=\ell_2=x$ and $h_1=h_2=y^2$.
Conditions $(M1)$ and $(M2)$
are verified as in the previous case, so this possibility occurs.\\ \\
$\bullet$ $G^\circ=\mZ_2\times\mZ_4, \; \;G=\mZ_2\times
Q_8$ and $G^\circ=\mZ_2\times\mZ_4, \; \;G=D_{4,4,-1}$. \\
All elements of order $2$ in $G$ are central so condition $(M1)$
cannot be satisfied and these cases do not occur.\\ \\
%
%
$\bullet$ $G^\circ=(\mZ_2)^3, \; \;G=\mZ_4\ltimes(\mZ_2)^2=\langle x,y,z \; | \; x^4=y^2=z^2=1,
xyx^{-1}=yz,\;[x,z]=[y,z]=1\rangle$\\
We consider the subgroup $G^\circ:=\langle
y,\,z,\,x^2\rangle\cong(\mZ_2)^3$. Set
$\ell_1=\ell_2=y$ and $h_1=z,\;\;h_2=x^2$.
Condition $(M1)$ holds because $G^\circ$ is abelian and $[x, y] \neq 1$. Condition $(M2)$
is satisfied because if $g\in xG^\circ$
then $g^2\in \langle z,\,x^2\rangle$. Therefore this case occurs.
\end{proof}

\subsection{The case $g(C)=7, \; |G|=36$} \label{g(C)=7}

\begin{proposition} \label{no-mixed-7}
The case $g(C)=7, \; |G|=36$ does not occur.
\end{proposition}
\begin{proof}
In this case $\mathbf{n}=(3)$, so $G^\circ$ is a group of order $18$
 which is $(1 \; |\; 3)$-generated. There are five groups of order $18$ up to
isomorphism. By computer search or direct calculation we see that
the only one which is $(1 \;| \; 3)$-generated is $\mZ_3 \times S_3=G(18,3)$. Thus
$G$ would fit into a short exact sequence
\begin{equation} \label{nonsplit7}
1 \lr \mZ_3 \times S_3 \lr G \lr \mZ_2 \lr 1.
\end{equation}
A computer search shows that the only groups of order $36$
containing a subgroup isomorphic to $\mZ_3 \times S_3$ are
$G(36,10)=S_3 \times S_3$  and $G(36,12)=\mZ_6 \times S_3$ (see
\verb|GAP4| script $2$ in the Appendix). They contain $15$ and $7$
elements of order $2$, respectively. On the other hand $\mZ_3 \times
S_3$ contains
 $3$ elements of order $2$, so by Remark \ref{elements-order-2}
 all possible extensions of the
form (\ref{nonsplit7}) are split and this case cannot occur.
\end{proof}

\subsection{The case $g(C)=9, \; |G|=64$} \label{g(C)=9}

\begin{proposition} \label{no-mixed-9}
The case $g(C)=9, \; |G|=64$ does not occur.
\end{proposition}
The proof will be the consequence of the results below. First notice that, since
$\mathbf{n}=(2)$, the group $G^\circ$ must be $(1\, | \,2)$-generated.
\begin{comp} \label{comp:1}
There exist precisely $8$ groups of order $32$ which are $(1 \; | \;
2)$-generated, namely $G(32,t)$ for $t \in \{2,4,5,6,7,8,12,17 \}$.
The number $n_2$ of their elements of order $2$ is given in the
table below:
\end{comp}
$$
\begin{array}{c|c|c|c|c|c|c|c|c}
 t &2&4&5&6&7&8&12&17\\
\hline n_2(G(32,t))&7&3&7&11&11&3&3&3\\
\end{array}
$$
\begin{proof}[Proof.]
Slightly modifying the first part of \verb|GAP4| script $1$ in the
Appendix we easily find that the groups of
order $32$ which are $(1\; | \;2)$-generated are exactly those in
the statement.
%
The number of elements of order $2$ in each case are found by a quick
computer search: see again the Appendix, \verb|GAP4| script $3$.
\end{proof}

\begin{comp} \label{comp:2}
Let $t \in \{2,4,5,6,7,8,12,17\}$. A nonsplit extension of the form
\begin{equation} \label{ext-nonsplit}
1 \lr G(32,t) \lr G(64,s) \lr \mZ_2 \lr 1
\end{equation}
exists if and only if the pair $(t,s)$ is one of the following: \\
\\
$(2,9)$, $(2,57)$, $(2,59)$, $(2,63)$, $(2,64)$, $(2,68)$, $(2,70)$,
$(2,72)$, $(2,76)$, $(2,79)$, $(2,81)$, $(2,82)$, \\
$(4,11)$, $(4,28)$, $(4,122)$, $(4,127)$, $(4,172)$, $(4,182)$, \\
$(5,5)$, $(5,9)$, $(5,112)$, $(5,113)$, $(5,114)$, $(5,132)$,
$(5,164)$, $(5,165)$, $(5,166)$, \\
$(6,33)$, $(6,35)$, \\
$(7,33)$, \\
$(8,37)$, \\
$(12,7)$, $(12,13)$, $(12,14)$, $(12,15)$, $(12,16)$, $(12,126)$,
$(12,127)$, $(12,143)$, $(12,156)$,\\ $(12,158)$, $(12,160)$, \\
$(17,28)$, $(17,43)$, $(17,45)$, $(17,46)$.
\end{comp}
\begin{proof}[Proof.]
Assume $t=2$. Using the \verb|GAP4| script $4$ in the Appendix we
find that the groups of order $64$ containing a subgroup isomorphic
to $G(32,2)$ are $G(64,s)$ for $s \in
\{8,\,9,\,56,\,57,\,58,\,59,\,61,\,62,\,63$, $64,\,
66,\,67,\,68,\,69,\,70,\,72, \,73,\,74,\,75,\,76,\,
77,\,78,\,79,\,80,\,81,\,82\}$. By Remark \ref{elements-order-2} and
Computational Fact \ref{comp:1}, in order to detect all the groups
$G(64,s)$ fitting in some nonsplit extension of type
\eqref{ext-nonsplit} with $t=2$, it is sufficient to select  from
the previous list the groups containing exactly $n_2=7$ elements of
order $2$. This can be done with the \verb|GAP4| script $5$ in the
Appendix, proving the claim in the case $t=2$. The proof
 for the other values of $t$ may be carried out exactly in the same way.
\end{proof}
Let us denote by  $[G,\,G]_2$ and $[G^\circ, \, G^\circ]_2$
the subsets of elements of order $2$ in $[G,\,G]$ and
$[G^\circ, G^\circ]$, respectively.
\begin{lemma} \label{deriv-2}
Assume $g(C)=9$ and that one of the following situations occur:
\begin{itemize}
\item  $[G,G]_2 \subseteq Z(G)$;
\item  there exists some element $y \in G \setminus G^\circ$
commuting with all elements in $[G^\circ, \, G^\circ]_2$.
\end{itemize}
Then given any generating vector $\mathcal{V}=\{\ell_1; \; h_1, h_2 \}$
of type $(1 \, | \, 2)$ for $G^\circ$,
condition $(M1)$ in
Proposition $\ref{building-mixed}$ cannot be satisfied.
\end{lemma}
\begin{proof}
Since $\ell_1 \in [G^\circ, \, G^\circ]_2 \subseteq
[G, \, G]_2$, in any of the above situations $C_G(\ell_1)$ is not
 contained in $G^\circ$, so $(M1)$ cannot hold.
\end{proof}

\begin{comp} \label{comp:3}
Let $G=G(64,s)$ be one of the groups appearing in the list of Computational Fact
 $\ref{comp:2}$. Then $[G,G]_2$ is not contained in $Z(G)$ if and only if
 $s=5$, $33$, $35$, $37$.
\end{comp}
\begin{proof}[Proof]
See the \verb|GAP4| script $6$ in the Appendix.
\end{proof}
Computational facts \ref{comp:2}, \ref{comp:3} and Lemma
\ref{deriv-2} imply that we only need to analyze the following pairs
$(G^\circ, \,G)$:
$$
\begin{array}{c|c}
G^\circ & G \\
\hline
G(32,5)& G(64,5) \\
G(32,6) & G(64,33) \\
G(32,7) & G(64,33) \\
G(32,6) & G(64,35) \\
G(32,8) & G(64,37)
\end{array}
$$
\begin{proposition} \label{no-32-5}
The case $G^\circ=G(32,5)$ does not occur.
\end{proposition}
\begin{proof}
A presentation for the group $G^\circ$ is
\begin{equation*}
G^\circ=\langle x,y,z~|~x^8=y^2=z^2=1, \;[y,z]=[x,\,z]=1,\;[x,\,y]=z\rangle.
\end{equation*}
Its derived subgroup contains exactly one element of
 order $2$, namely $z$. It follows that if $\{\ell_1; \, h_1, h_2 \}$ is any generating
  vector of type $(1 \, | \, 2)$ for $G^{\circ}$, then $\ell_1=z$. Since $[G^\circ, \, G^\circ]$
 is characteristic in $G^\circ$, condition $(M1)$ cannot be satisfied for any
 embedding of $G^\circ$ into $G$.
\end{proof}
By using the two instructions
\verb|P:=PresentationViaCosetTable(G)| and \verb|TzPrintRelators(P)|
and setting in the output
\begin{equation*}
x:=\verb|f1|, \; y:=\verb|f2|, \; z:=\verb|f3|, \; w:=\verb|f4|, \;
v:=\verb|f5|, \; u:=\verb|f6|
\end{equation*}
one obtains the following presentations for $G(64,33)$, $G(64,35)$ and
$G(64,37)$.
\begin{equation} \label{G(64,33)}
\begin{split}
G(64,33)=\langle
  x,\,y,\,z,\,w,\,v,\,u~|&~z^2=w^2=v^2=u^2=1,\;x^2=w,\;y^2=u,\\
&\;[x,\,zy]=z,\;[x,\,vz]=v,\;[x,\,vu]=u,\\
&\;[y,\,z]=[y,\,v]=[z,\,v]=[w,\,v]=[x,\,u]=1\rangle
\end{split}
\end{equation}

\begin{equation} \label{G(64,35)}
\begin{split}
G(64,35)=\langle
  x,\,y,\,z,\,w,\,v,\,u~|&~w^2=v^2=u^2=1,\; z^2=y^2=u,\;x^2=w,\\
&\;[y,\,z]=[z,\,w]=u,\,[x, yz]=z,\,[x,\,z]=uv,\\
&\;[y,\,v]=[z,\,v]=[w,\,v]=[x,\,u]=1\rangle
\end{split}
\end{equation}

\begin{equation} \label{G(64,37)}
\begin{split}
G(64,37)=\langle
  x,\,y,\,z,\,w,\,v,\,u~|&~v^2=u^2=1,\;w^2=z^2=y^2=u,\;x^2=w,\\
&\;[y,\,z]=[z,\,w]=u,\,[x,\,yz]=z,\,[x,\,z]=uv,\\
&\;[y,\,v]=[z,\,v]=[w,\,v]=1\rangle.
\end{split}
\end{equation}

\begin{comp} \label{comp:4}
Referring to presentations \eqref{G(64,33)}, \eqref{G(64,35)} and
 \eqref{G(64,37)}, we have the following facts.
\begin{itemize}
\item The group $G(64,33)$ contains exactly one subgroup $N_1$ isomorphic to $G(32,6)$ and one
 subgroup $N_2$ isomorphic to $G(32,7)$, namely
\begin{equation*}
N_1:=\langle x,\,z,\,w,\,v,\,u \rangle, \quad N_2:=\langle xy,\,z,\,w,\,v,\,u \rangle.
\end{equation*}
\item The group $G(64,35)$ contains exactly two subgroups $N_3$, $N_4$ isomorphic
to $G(32,6)$, namely
\begin{equation*}
N_3:=\langle x,\,z,\,w,\,v,\,u \rangle, \quad N_4:=\langle xy,\,z,\,w,\,v,\,u \rangle.
\end{equation*}
\item The group $G(64,37)$ contains exactly two subgroups $N_5$, $N_6$ isomorphic
to $G(32,8)$, namely
\begin{equation*}
N_5:=\langle x,\,z,\,w,\,v,\,u \rangle, \quad N_6:=\langle xy,\,z,\,w,\,v,\,u \rangle.
\end{equation*}
\end{itemize}
In addition, for every $i \in \{1, \ldots, 6 \}$ we have
\begin{itemize}
\item[$(a)$]
$[N_i, \, N_i]= \langle v,u \rangle \cong \mZ_2 \times \mZ_2$.
\item[$(b)$]
$y \notin N_i$ and $y$ commutes with all elements in $[N_i, \, N_i]$.
\end{itemize}
\end{comp}
\begin{proof}[Proof]
See the \verb|GAP4| script $7$ in the Appendix.
\end{proof}

\begin{proposition} \label{no-32-6-7-8}
The cases $G^\circ=G(32,6), \, G(32,7), \, G(32,8)$ do not occur.
\end{proposition}
\begin{proof}
By Lemma \ref{deriv-2} and Computational Fact \ref{comp:4} it follows that,
 given any nonsplit extension of type \eqref{ext-nonsplit} with $G^\circ$ as above,
 condition $(M1)$ in Proposition \ref{building-mixed} cannot be satisfied.
 \end{proof}
Summing up, we finally obtain \\ \\
\emph{Proof of Proposition }\ref{no-mixed-9}. It follows from Propositions \ref{no-32-5} and
 \ref{no-32-6-7-8}.
\begin{flushright}
$\Box$
\end{flushright}

\section{Moduli spaces} \label{moduli spaces}
Let $\mathfrak{M}_{a,b}$ be the moduli space of smooth minimal
surfaces of general type with $\chi(\mO_S)=a, \; K_S^2=b$; by an
important result of Gieseker, $\mathfrak{M}_{a,b}$ is a
quasiprojective variety for all $a,b \in \mathbb{N}$ (see
\cite{Gie77}). Obviously, our surfaces are contained in
$\mathfrak{M}_{1,8}$ and we want to describe their locus there. We
denote by $\mathfrak{M}(G, \mathbf{m}, \mathbf{n})$ the moduli space
of unmixed surfaces of type $(G, \mathbf{m}, \mathbf{n})$ and by
$\mathfrak{M}(G, \mathbf{n})$ the moduli space of mixed surfaces
 of type $(G, \mathbf{n})$. We know that
  $\mathbf{n}=(2^2)$, $(3)$ or $(2)$
 in the unmixed case, whereas $\mathbf{n}=(2^2)$ in the mixed one.
 By a general result of Catanese (\cite{Ca00}), both
  $\mathfrak{M}(G, \mathbf{m}, \mathbf{n})$ and $\mathfrak{M}(G, \mathbf{n})$
 consist of finitely many irreducible connected
components of $\mathfrak{M}_{1,8}$, all of the same dimension. More
precisely, we have
\begin{equation*}
\textrm{dim \,} \mathfrak{M}(G, \mathbf{m}, \mathbf{n})=r+s-3, \quad
\textrm{dim \,} \mathfrak{M}(G,\mathbf{n})=s.
\end{equation*}
Consider the mapping class groups in genus zero and one:
\begin{equation*}
\begin{split}
\textrm{Mod}_{0,[r]}:= \langle \sigma_1, \ldots, \sigma_r \; | \;&
\sigma_i \sigma_{i+1} \sigma_i=\sigma_{i+1} \sigma_i
\sigma_{i+1},\\
& \sigma_i \sigma_j = \sigma_j \sigma_i \quad \textrm{if} \; |i-j| \geq 2,\\
& \sigma_{r-1} \sigma_{r-2} \cdots \sigma_1^2 \cdots \sigma_{r-2}
\sigma_{r-1}=1 \rangle,
\end{split}
\end{equation*}
\begin{equation*}
\textrm{Mod}_{1,1}:= \langle t_{\alpha}, t_{\beta}, t_{\gamma} \; |
\;  t_{\alpha}t_{\beta}t_{\alpha}=t_{\beta}t_{\alpha}t_{\beta}, \;
\; (t_{\alpha}t_{\beta})^3=1 \rangle,
\end{equation*}
\begin{equation*}
\begin{split}
\textrm{Mod}_{1,[2]}:= \langle t_{\alpha}, t_{\beta}, t_{\gamma},
\rho \; | \;
t_{\alpha}t_{\beta}t_{\alpha}&=t_{\beta}t_{\alpha}t_{\beta}, \; \;
t_{\alpha}t_{\gamma}t_{\alpha}=t_{\gamma}t_{\alpha}t_{\gamma}, \\
t_{\beta}t_{\gamma}&=t_{\gamma}t_{\beta}, \; \;
(t_{\alpha}t_{\beta}t_{\gamma})^4=1,\\
t_{\alpha}\rho&=\rho t_{\alpha}, \; \; t_{\beta}\rho=\rho t_{\beta},
\; \; t_{\gamma}\rho=\rho t_{\gamma} \rangle.
\end{split}
\end{equation*}
One can prove that
\begin{equation*}
\begin{split}
\textrm{Mod}_{0,[r]}:&=\pi_0\; \textrm{Diff}^+(\mathbb{P}^1-\{p_1, \ldots, p_r \}), \\
\textrm{Mod}_{1,1}:&=\pi_0\; \textrm{Diff}^+(\Sigma_1-\{p \}), \\
\textrm{Mod}_{1,[2]}:&=\pi_0\; \textrm{Diff}^+(\Sigma_1-\{p,\,q \}),
\end{split}
\end{equation*}
where $\Sigma_1$ is the torus $S^1 \times S^1$ (\cite{Schn03},
\cite{CattMu04}). This implies that we can define actions of these
groups on the set of generating vectors for $G$ of type $(0 \; | \;
m_1, \ldots, m_r)$, $(1 \; | \; n)$ and
$(1 \; | \; n^2)$, respectively. \\ \\
If $\mathcal{V}:=\{g_1, \ldots, g_r\}$ is of type $(0 \; | \; m_1,
\ldots, m_r)$ then the action is given by
\begin{equation*}
\sigma_i: \left \{
\begin{array}{ll}
g_i & \lr g_{i+1} \\
g_{i+1} & \lr g_{i+1}^{-1} g_i g_{i+1} \\
g_j & \lr g_j  \quad \textrm{if} \; j \neq i,\; i+1.
\end{array} \right.
\end{equation*}
If $\mathcal{W}:=\{\ell_1; \, h_1, h_2\}$ is of type $(1 \; | \; n)$
 then
\begin{equation*}
t_{\alpha} \colon \left\{ \begin{array}{ll}
\ell_1 \lr \ell_1 \\
h_1 \lr h_1 \\
h_2 \lr h_2h_1 \\
\end{array} \right.
\quad t_{\beta} \colon \left\{ \begin{array}{ll}
\ell_1 \lr \ell_1 \\
h_1 \lr h_1h_2^{-1} \\
h_2 \lr h_2. \\
\end{array} \right.
\end{equation*}
If $\mathcal{W}:=\{\ell_1, \ell_2; \, h_1, h_2\}$ is of type $(1 \;
| \; n^2)$
 then
\begin{equation*}
\begin{split}
& t_{\alpha} \colon \left\{ \begin{array}{ll}
  \ell_1 \lr \ell_1 \\
  \ell_2 \lr \ell_2 \\
  h_1 \lr h_1 \\
  h_2 \lr h_2h_1
\end{array} \right.
\quad  \quad \quad \quad \quad \quad \quad t_{\beta} \colon \left\{
\begin{array}{ll}
  \ell_1 \lr \ell_1 \\
  \ell_2 \lr \ell_2\\
  h_1 \lr h_1h_2^{-1} \\
  h_2 \lr h_2
\end{array} \right.
\\
& t_{\gamma} \colon \left\{ \begin{array}{ll}
  \ell_1 \lr \ell_1 \\
  \ell_2 \lr h_1h_2^{-1}h_1^{-1}\ell_2h_1h_2h_1^{-1}\\
  h_1 \lr h_2^{-1}\ell_1h_1 \\
  h_2 \lr h_2
\end{array} \right.
\quad  \; \rho \colon \left\{ \begin{array}{ll}
  \ell_1 \lr h_2^{-1}h_1^{-1}\ell_2 h_1h_2 \\
  \ell_2 \lr h_1^{-1}h_2^{-1}\ell_1 h_2h_1 \\
  h_1 \lr h_1^{-1} \\
  h_2 \lr h_2^{-1}.
\end{array} \right.
\end{split}
\end{equation*}
These are called \emph{Hurwitz moves} and the induced equivalence
relation on generating vectors
is said \emph{Hurwitz equivalence} (see \cite{BaCa03},
\cite{BaCaGr06}, \cite{Pol07}). \\ \\
Now let $\mathfrak{B}(G, \, \mathbf{m}, \, \mathbf{n})$ be the set
of pairs of generating vectors $(\mathcal{V}, \; \mathcal{W})$ such
that the assumptions in Proposition \ref{structureresult2} are
satisfied; then we denote by $\mathfrak{R}$ the equivalence relation
on $\mathfrak{B}(G, \, \mathbf{m}, \, \mathbf{n})$ generated by
Hurwitz moves on $\mathcal{V}$, Hurwitz moves on $\mathcal{W}$ and the
simultaneous action of $\textrm{Aut}(G)$ on $\mathcal{V}$ and $\mathcal{W}$.
Similarly, let $\mathfrak{B}(G, \,\mathbf{n})$ be the set of
generating vectors $\mathcal{V}$ such that the assumptions of Proposition
\ref{building-mixed} are satisfied; then we denote by $\mathfrak{R}$
the equivalence relation on $\mathfrak{B}(G, \,\mathbf{n})$ generated
by the Hurwitz moves and the action of $\textrm{Aut}(G)$ on $\mathcal{V}$.
\begin{proposition} \label{prop:moduli spaces}
The number of irreducible components in $\mathfrak{M}(G, \,
\mathbf{m}, \, \mathbf{n})$ equals the number of
$\mathfrak{R}$-classes in $\mathfrak{B}(G, \, \mathbf{m}, \, \mathbf{n})$.
Analogously, the number of irreducible components in $\mathfrak{M}(G, \,
\mathbf{n})$ equals the number of $\mathfrak{R}$-classes in
$\mathfrak{B}(G, \, \mathbf{n})$.
\end{proposition}
\begin{proof}
We can repeat exactly the same argument used in \cite[Propositions
5.2 and 5.5]{BaCaGr06}; we must just replace, where it is necessary,
the mapping class group of $\mathbb{P}^1$ with the mapping class
group of the elliptic curve $E$.
\end{proof}
Proposition \ref{prop:moduli spaces} in principle allows us to
compute the number of connected components of the moduli space in
each case. In practice, this task may be too hard to be achieved by
hand, but it is not out of reach if one uses the
computer. Recently, M. Penegini and S. Rollenske developed a
\verb|GAP4| script that solves this problem in a rather short time.
We put the result of their calculations in the Main Theorem (see
Introduction), referring the reader to the forthcoming
 paper \cite{Pe08} for further details.

\section{Appendix}

In this Appendix we include, for the reader's convenience, some of
the \verb|GAP4| scripts that we have used in our computations; all
the others
are similar and can be easily obtained modifying the ones below. \\
\\ Let us show how the procedure in the proof of Proposition
\ref{g(F)=3nonmixed} applies to an explicit example, namely
$\mathbf{m}_{\alpha(\mathbf{m})}=(2,4,12)_{12}$.  First we find all the nonabelian groups of order
$24$ that are $(0 \, | \,2,4,12)$-generated. This is done using \verb|GAP4| as below;
 the output tells us that there is only one
such a group, namely $G=G(24,5)$. \bigskip
\begin{verbatim}
gap> # -------------- SCRIPT 1 ------------------
gap> s:=NumberSmallGroups(24);; set:=[1..s];
[1..15]
gap> for t in set do
> c:=0; G:= SmallGroup(24,t);
> Ab:=IsAbelian(G);
> for g1 in G do
> for g2 in G do
> g3:=(g1*g2)^-1;
> H:= Subgroup(G, [g1,g2]);
> if Order(g1)=2 and Order(g2)=4 and Order(g3)=12 and
> Order(H)=Order(G) and
> Ab=false then
> c:=c+1; fi;
> if Order(g1)=2 and Order(g2)=4 and Order(g3)=12 and
> Order(H)=Order(G) and
> Ab=false and c=1 then
> Print(IdSmallGroup(G)," ");
> fi; od; od; od; Print("\n");
[24,5]
\end{verbatim}
\bigskip
By using the two instructions \verb|P:=PresentationViaCosetTable(G)|
and \verb|TzPrintRelators(P)| we see that $G$ has the presentation
 $\langle x, y \; | \; x^2=y^{12}=1, \; xyx^{-1}=y^5
\rangle$, hence it is isomorphic to the metacyclic group
$D_{2,12,5}$. \\In order to speed up further computations, we
define the sets $G2$, $G4$ given by the elements of $G$ having
order $2$ and $4$, respectively.
\bigskip
\begin{verbatim}
gap> G:=SmallGroup(24,5);;
gap> G2:=[];; G4:=[];;
gap> for g in G do
> if Order(g)=2 then Add(G2,g); fi;
> if Order(g)=4 then Add(G4,g); fi; od;
\end{verbatim}
\bigskip
Then we check whether
 $G$ is actually $(1 \, | \,2^2)$-generated; if not, it should
 be excluded.
\bigskip
\begin{verbatim}
gap> c:=0;;
gap> for l2 in G2 do
> for h1 in G do
> for h2 in G do
> l1:=(l2*h1*h2*h1^-1*h2^-1)^-1;
> K:=Subgroup(G, [l2, h1, h2]);
> if Order(l1)=2 and Order(K)=Order(G) then
> Print(IdSmallGroup(G), " is (1 | 2,2)-generated", "\n"); c:=1; fi;
> if c=1 then break; fi; od;
> if c=1 then break; fi; od;
> if c=1 then break; fi; od;
[24,5] is (1 | 2,2)-generated
\end{verbatim}
\bigskip
We finish the proof by checking whether the surface $S$ actually exists;
the procedure is to look for a pair $(\mathcal{V}, \, \mathcal{W})$ of
generating vectors for $G$ satisfying the assumptions of Proposition
\ref{structureresult2}.
\begin{verbatim}
gap> c:=0;;
gap> for g1 in G2 do
> for g2 in G4 do
> g3:=(g1*g2)^-1;
> H:=Subgroup(G, [g1, g2]);
> for l2 in G2 do
> for h1 in G do
> for h2 in G do
> l1:=(l2*h1*h2*h1^-1*h2^-1)^-1;
> K:=Subgroup(G, [l2, h1, h2]);
> Boole1:=l1 in ConjugacyClass(G, g1);
> Boole2:=l1 in ConjugacyClass(G, g2^2);
> Boole3:=l1 in ConjugacyClass(G, g3^6);
> Boole4:=l2 in ConjugacyClass(G, g1);
> Boole5:=l2 in ConjugacyClass(G, g2^2);
> Boole6:=l2 in ConjugacyClass(G, g3^6);
> if Order(g3)=12 and Order(l1)=2 and
> Order(H)=Order(G) and Order(K)=Order(G) and
> Boole1=false and Boole2=false and Boole3=false and
> Boole4=false and Boole5=false and Boole6=false then
> Print("The surface exists "); c:=1; fi;
> if c=1 then break; fi; od;
> if c=1 then break; fi; od;
> if c=1 then break; fi; od;
> if c=1 then break; fi; od;
> if c=1 then break; fi; od; Print("\n");
The surface exists
\end{verbatim}
\bigskip
The script above can be easily modified in order to obtain the list
of all admissible pairs $(\mathcal{V}, \,\mathcal{W})$; for
instance, one of such pairs is given by
\begin{equation*}
\begin{split}
g_1&=x, \; \; g_2=xy^{-1}, \;\; g_3=y \\
\ell_1&=xy^2, \; \; \ell_2=xy^2, \; \; h_1=y, \;\; h_2=y.
\end{split}
\end{equation*}

 \bigskip \bigskip

Finally, here are the \verb|GAP4| scripts used in Section \ref{classification-mix}.
\bigskip \medskip
\begin{verbatim}
gap> # -------------- SCRIPT 2 ------------------
gap> s:=NumberSmallGroups(36);; set:=[1..s];
[1..14]
gap> for t in set do
> c:=0; G:=SmallGroup(36,t);
> N:=NormalSubgroups(G);
> for G0 in N do
> if IdSmallGroup(G0)=[18,3] then
> c:=c+1; fi;
> if IdSmallGroup(G0)=[18,3] and c=1 then
> Print(IdSmallGroup(G), " ");
> fi; od; od; Print("\n");
[36,10] [36,12]
\end{verbatim}

\bigskip \medskip

\begin{verbatim}
gap> # -------------- SCRIPT 3 ------------------
gap> set:=[2,4,5,6,7,8,12,17];;
gap> for t in set do
> n2:=0;
> G0:=SmallGroup(32,t);
> for g in G0 do
> if Order(g)=2 then
> n2:=n2+1; fi; od;
> Print(IdSmallGroup(G0), " "); Print(n2, "    ");
> od; Print("\n");
[32,2] 7    [32,4] 3    [32,5] 7    [32,6] 11    [32,7] 11
[32,8] 3    [32,12] 3    [32,17] 3
\end{verbatim}

\bigskip \medskip

\begin{verbatim}
gap> # -------------- SCRIPT 4 ------------------
gap> s:=NumberSmallGroups(64);; set:=[1..s];
[1..267]
gap> for t in set do
> c:=0; G:=SmallGroup(64,t);
> N:=NormalSubgroups(G);
> for G0 in N do
> if IdSmallGroup(G0)=[32,2] then
> c:=c+1; fi;
> if IdSmallGroup(G0)=[32,2] and c=1 then
> Print(IdSmallGroup(G), " ");
> fi; od; od; Print("\n");
[64,8] [64,9] [64,56] [64,57] [64,58] [64,59] [64,61] [64,62] [64,63] [64,64]
[64,66] [64,67] [64,68] [64,69] [64,70] [64,72] [64,73] [64,74] [64,75] [64,76]
[64,77] [64,78] [64,79] [64,80] [64,81] [64,82]
\end{verbatim}

\bigskip \medskip


\begin{verbatim}
gap> # -------------- SCRIPT 5 ------------------
gap> set:=[8,9,56,57,58,59,61,62,63,64,66,67,68,69,70,
>72,73,74,75,76,77,78,79,80,81,82];;
gap> for t in set do
> n2:=0; G:=SmallGroup(64,t);
> for g in G do
> if Order(g)=2 then n2:=n2+1;
> fi; od;
> if n2=7 then
> Print(IdSmallGroup(G), " ");
> fi; od; Print("\n");
[64,9] [64,57] [64,59] [64,63] [64,64] [64,68] [64,70]
[64,72] [64,76] [64,79] [64,81] [64,82]
\end{verbatim}

\bigskip \medskip


\begin{verbatim}
gap> # -------------- SCRIPT 6 ------------------
gap> set:=[5,7,9,11,13,,14,15,16,28,33,35,37,43,45,46,
>57,59,63,64,68,70,72,76,79,81,82,112,113,114, 122,126,
>127,132,143,156,158,160,164,165,166,172,182];;
gap> for t in set do
> c:=0; G:=SmallGroup(64,t);
> D:=DerivedSubgroup(G);
> for d in D do
> B:=d in Center(G);
> if Order(d)=2 and B=false then
> c:=c+1; fi;
> if Order(d)=2 and B=false and c=1 then
> Print(IdSmallGroup(G), " ");
> fi; od; od; Print("\n");
[64,5] [64,33] [64,35] [64,37]
\end{verbatim}

\bigskip \medskip


\begin{verbatim}
gap> # -------------- SCRIPT 7 ------------------
gap> s:=[33, 35, 37];; I:=[1, 2, 3];;
gap> r:=[ [[32,6], [32,7]], [[32,6]], [[32,8]] ];;
> for i in I do
> G:=SmallGroup(64, s[i]); Print(IdSmallGroup(G), "\n");
> for N in NormalSubgroups(G) do
> if IdSmallGroup(N) in r[i] then
> Print(N, "="); Print(IdSmallGroup(N), "   ");
> Print(DerivedSubgroup(N), "\n");
> fi; od; Print("\n"); od;
[64,33]
Group( [ f1*f2, f3, f4, f5, f6 ] )=[32,7]   Group( [ f5, f6 ] )
Group( [ f1, f3, f4, f5, f6 ] )=[32,6]   Group( [ f5, f6 ] )

[64,35]
Group( [ f1*f2, f3, f4, f5, f6 ] )=[32,6]   Group( [ f5, f6 ] )
Group( [ f1, f3, f4, f5, f6 ] )=[32,6]   Group( [ f5, f6 ] )

[64,37]
Group( [ f1*f2, f3, f4, f5, f6 ] )=[32,8]   Group( [ f5, f6 ] )
Group( [ f1, f3, f4, f5, f6 ] )=[32,8]   Group( [ f5, f6 ] )
\end{verbatim}

\end{document}